\documentstyle[11pt,amstex,amscd,graphicx]{amsart}

\begin{document}

\newtheorem{theorem}{Theorem}[section]
\newtheorem{prop}[theorem]{Proposition}
\newtheorem{lemma}[theorem]{Lemma}
\newtheorem{cor}[theorem]{Corollary}
\newtheorem{definition}[theorem]{Definition}
\newtheorem{conj}[theorem]{Conjecture}
\newtheorem{claim}[theorem]{Claim}

\newcommand{\boundary}{\partial}
\newcommand{\C}{{\mathbb C}}
\newcommand{\integers}{{\mathbb Z}}
\newcommand{\natls}{{\mathbb N}}
\newcommand{\ratls}{{\mathbb Q}}
\newcommand{\reals}{{\mathbb R}}
\newcommand{\proj}{{\mathbb P}}
\newcommand{\lhp}{{\mathbb L}}
\newcommand{\tube}{{\mathbb T}}
\newcommand{\cusp}{{\mathbb P}}
\newcommand\AAA{{\mathcal A}}
\newcommand\BB{{\mathcal B}}
\newcommand\CC{{\mathcal C}}
\newcommand\DD{{\mathcal D}}
\newcommand\EE{{\mathcal E}}
\newcommand\FF{{\mathcal F}}
\newcommand\GG{{\mathcal G}}
\newcommand\HH{{\mathcal H}}
\newcommand\II{{\mathcal I}}
\newcommand\JJ{{\mathcal J}}
\newcommand\KK{{\mathcal K}}
\newcommand\LL{{\mathcal L}}
\newcommand\MM{{\mathcal M}}
\newcommand\NN{{\mathcal N}}
\newcommand\OO{{\mathcal O}}
\newcommand\PP{{\mathcal P}}
\newcommand\QQ{{\mathcal Q}}
\newcommand\RR{{\mathcal R}}
\newcommand\SSS{{\mathcal S}}
\newcommand\TT{{\mathcal T}}
\newcommand\UU{{\mathcal U}}
\newcommand\VV{{\mathcal V}}
\newcommand\WW{{\mathcal W}}
\newcommand\XX{{\mathcal X}}
\newcommand\YY{{\mathcal Y}}
\newcommand\ZZ{{\mathcal Z}}
\newcommand\CH{{\CC\HH}}
\newcommand\MF{{\MM\FF}}
\newcommand\PMF{{\PP\kern-2pt\MM\FF}}
\newcommand\ML{{\MM\LL}}
\newcommand\PML{{\PP\kern-2pt\MM\LL}}
\newcommand\GL{{\GG\LL}}
\newcommand\Pol{{\mathcal P}}
\newcommand\half{{\textstyle{\frac12}}}
\newcommand\Half{{\frac12}}
\newcommand\Mod{\operatorname{Mod}}
\newcommand\Area{\operatorname{Area}}
\newcommand\ep{\epsilon}
\newcommand\hhat{\widehat}
\newcommand\Proj{{\mathbf P}}
\newcommand\U{{\mathbf U}}
 \newcommand\Hyp{{\mathbf H}}
\newcommand\D{{\mathbf D}}
\newcommand\Z{{\mathbb Z}}
\newcommand\R{{\mathbb R}}
\newcommand\Q{{\mathbb Q}}
\newcommand\E{{\mathbb E}}
\newcommand\til{\widetilde}
\newcommand\length{\operatorname{length}}
\newcommand\tr{\operatorname{tr}}
\newcommand\gesim{\succ}
\newcommand\lesim{\prec}
\newcommand\simle{\lesim}
\newcommand\simge{\gesim}
\newcommand{\simmult}{\asymp}
\newcommand{\simadd}{\mathrel{\overset{\text{\tiny $+$}}{\sim}}}
\newcommand{\ssm}{\setminus}
\newcommand{\diam}{\operatorname{diam}}
\newcommand{\pair}[1]{\langle #1\rangle}
\newcommand{\T}{{\mathbf T}}
\newcommand{\inj}{\operatorname{inj}}
\newcommand{\pleat}{\operatorname{\mathbf{pleat}}}
\newcommand{\short}{\operatorname{\mathbf{short}}}
\newcommand{\vertices}{\operatorname{vert}}
\newcommand{\collar}{\operatorname{\mathbf{collar}}}
\newcommand{\bcollar}{\operatorname{\overline{\mathbf{collar}}}}
\newcommand{\I}{{\mathbf I}}
\newcommand{\tprec}{\prec_t}
\newcommand{\fprec}{\prec_f}
\newcommand{\bprec}{\prec_b}
\newcommand{\pprec}{\prec_p}
\newcommand{\ppreceq}{\preceq_p}
\newcommand{\sprec}{\prec_s}
\newcommand{\cpreceq}{\preceq_c}
\newcommand{\cprec}{\prec_c}
\newcommand{\topprec}{\prec_{\rm top}}
\newcommand{\Topprec}{\prec_{\rm TOP}}
\newcommand{\fsub}{\mathrel{\scriptstyle\searrow}}
\newcommand{\bsub}{\mathrel{\scriptstyle\swarrow}}
\newcommand{\fsubd}{\mathrel{{\scriptstyle\searrow}\kern-1ex^d\kern0.5ex}}
\newcommand{\bsubd}{\mathrel{{\scriptstyle\swarrow}\kern-1.6ex^d\kern0.8ex}}
\newcommand{\fsubeq}{\mathrel{\raise-.7ex\hbox{$\overset{\searrow}{=}$}}}
\newcommand{\bsubeq}{\mathrel{\raise-.7ex\hbox{$\overset{\swarrow}{=}$}}}
\newcommand{\tw}{\operatorname{tw}}
\newcommand{\base}{\operatorname{base}}
\newcommand{\trans}{\operatorname{trans}}
\newcommand{\rest}{|_}
\newcommand{\bbar}{\overline}
\newcommand{\UML}{\operatorname{\UU\MM\LL}}
\newcommand{\EL}{\mathcal{EL}}
\newcommand{\tsum}{\sideset{}{'}\sum}
\newcommand{\tsh}[1]{\left\{\kern-.9ex\left\{#1\right\}\kern-.9ex\right\}}
\newcommand{\Tsh}[2]{\tsh{#2}_{#1}}
\newcommand{\qeq}{\mathrel{\approx}}
\newcommand{\Qeq}[1]{\mathrel{\approx_{#1}}}
\newcommand{\qle}{\lesssim}
\newcommand{\Qle}[1]{\mathrel{\lesssim_{#1}}}
\newcommand{\simp}{\operatorname{simp}}
\newcommand{\vsucc}{\operatorname{succ}}
\newcommand{\vpred}{\operatorname{pred}}
\newcommand\fhalf[1]{\overrightarrow {#1}}
\newcommand\bhalf[1]{\overleftarrow {#1}}
\newcommand\sleft{_{\text{left}}}
\newcommand\sright{_{\text{right}}}
\newcommand\sbtop{_{\text{top}}}
\newcommand\sbot{_{\text{bot}}}
\newcommand\sll{_{\mathbf l}}
\newcommand\srr{_{\mathbf r}}
\newcommand\geod{\operatorname{\mathbf g}}
\newcommand\mtorus[1]{\boundary U(#1)}
\newcommand\A{\mathbf A}
\newcommand\Aleft[1]{\A\sleft(#1)}
\newcommand\Aright[1]{\A\sright(#1)}
\newcommand\Atop[1]{\A\sbtop(#1)}
\newcommand\Abot[1]{\A\sbot(#1)}
\newcommand\boundvert{{\boundary_{||}}}
\newcommand\storus[1]{U(#1)}
\newcommand\Momega{\omega_M}
\newcommand\nomega{\omega_\nu}
\newcommand\twist{\operatorname{tw}}
\newcommand\modl{M_\nu}
\newcommand\MT{{\mathbb T}}
\newcommand\Teich{{\mathcal T}}
\renewcommand{\Re}{\operatorname{Re}}
\renewcommand{\Im}{\operatorname{Im}}

\title{A Combination Theorem for Strong Relative Hyperbolicity}

\author{Mahan Mj}
\address{RKM Vivekananda University, Belur Math, WB-711 202, India}

\author{Lawrence Reeves}
\address{University of Melbourne, Victoria 3010, Australia}
\date{}

\begin{abstract} 
We prove a combination theorem for trees of (strongly) relatively
hyperbolic spaces and finite graphs of (strongly) relatively
hyperbolic groups. This gives a geometric extension of Bestvina
and Feighn's Combination Theorem for hyperbolic groups and answers a question of Swarup. We also prove a converse
to the main Combination Theorem.
\end{abstract}

\maketitle

\begin{center}
AMS subject classification =   20F32(Primary), 57M50(Secondary)
\end{center}

\tableofcontents

\section{Introduction}
In \cite{BF}, Bestvina and Feighn proved a combination theorem for
hyperbolic groups. Motivated by this, Swarup asked the analogous
question \cite{bestvinahp} for relatively hyperbolic groups. Dahmani
\cite{dahmani-comb} and Alibegovic \cite{alib-comb} have proven
combination theorems motivated by applications to convergence groups
and
 limit groups
(cf. Sela \cite{sela-ltgps}). 

In this paper, we prove a {\bf geometric combination theorem} (as
opposed to a dynamical one) for
trees of
(strong) relatively hyperbolic metric spaces. We use Bestvina and
Feighn's Combination Theorem \cite{BF} directly in deducing the
relevant combination Theorem. The conditions we impose are
quite different from  those of \cite{dahmani-comb} and
\cite{alib-comb}. Our main Theorems \ref{strongcombin} 
and \ref{converse} are stated below:
 \medskip \medskip \\
 {\bf Strong Combination Theorem and converse: Theorems \ref{strongcombin} , \ref{converse}}
Let $X$ be a tree ($T$) of strongly relatively hyperbolic spaces
satisfying 
\begin{enumerate}
\item the q(uasi)-i(sometrically)-embedded condition 
\item  the strictly type-preserving
condition
\item the qi-preserving electrocution condition
\item the induced tree of coned-off spaces
satisfies the {\bf hallways flare} condition
\item the {\bf cone-bounded hallways strictly flare} condition.
\end{enumerate}
Then $X$  is {\em
  strongly hyperbolic} relative to the family $\CC$ of {\em maximal
  cone-subtrees of horosphere-like spaces}. 
\\ Conversely, if
$X$ be a tree ($T$) of strongly relatively hyperbolic spaces
satisfying conditions (1), (2), (3) such that 
$X$  is {\em
  strongly hyperbolic} relative to the family $\CC$ of {\em maximal
  cone-subtrees of horosphere-like spaces}, then the tree 
of spaces satisfies conditions (4), (5).

\medskip \medskip

Of the conditions given in the above theorem, Condition (1) is taken
directly from \cite{BF}. Condition (2) roughly says that the
 pre-image of a
horosphere-like subset (thought of as parabolic)
 in a vertex space (under the edge-space to
vertex-space map) is either empty or a 
horosphere-like subset in the corresponding edge-space. This condition
 may
be likened to the restriction to strictly type-preserving maps in the
theory of Kleinian groups. Condition (2) ensures an induced tree of
electrocuted spaces. Condition (3) says that the induced tree of
spaces also satisfies the qi-embedded condition. Condition (4) is
again taken directly from \cite{BF}. Condition (5) is the one
essential new condition. It says roughly that a pair of
geodesics whose vertices consist only of cone-points cannot lie close
to each other for long. 
The notion of {\em fully quasiconvex
  subgroups} introduced by Dahmani \cite{dahmani-comb}
 is related to Condition (3),
{\em the qi-preserving electrocution condition}. \\

\smallskip

\noindent {\bf Note:} {\it  In this paper
we adopt the convention that the  are 
horosphere-like subsets are {\bf coarsely proper}, i.e. no finite neighborhood of a horosphere-like subset is the whole space
(cf. \cite{bdm} which follows
a similar convention). This excludes the trivial case that $X$ is strongly
hyperbolic relative to itself, or a net in $X$. This assumption
translates
into the context of groups. Hence we assume that if a group
$G$ is strongly hyperbolic relative to a collection of
subgroups, then no subgroup $H$ in the collection is of finite
index in $G$.}

\smallskip

As an immediate consequence of theorem \ref{strongcombin}, we have: 
\\ \medskip \medskip
 
{\bf Strong Combination Theorem for Graphs of Groups: Theorem
  \ref{strongcombingp}} 
Let $G$ be a finite graph ($\Gamma$) of strongly relatively hyperbolic groups
satisfying 
\begin{enumerate}
\item the qi-embedded condition 
\item  the strictly type-preserving
condition
\item the qi-preserving electrocution condition
\item the induced tree of coned-off spaces
satisfies the {\bf hallways flare} condition
\item the {\bf cone-bounded hallways strictly flare} condition
\end{enumerate}
Then $G$  is {\em
  strongly hyperbolic} relative to the family $\CC$ of {\em maximal
  parabolic subgroups}.  \medskip

All these conditions are satisfied in the classical case of a
3-manifold fibering over the circle with fiber a punctured surface.
The one condition that needs checking is the {\em hallways flare
  condition} for the induced tree (in fact line) of coned-off spaces.
This fact is due to  Bowditch (see \cite{bowditch-ct} Section 6). The
verification involves using the associated singular structure coming
from stable and unstable foliations. We shall give a slightly modified
version, using an idea of Mosher \cite{mosher-hbh} to show this. (See
Section 4.3). In fact we prove the stronger Theorem: \\

\medskip

{\bf Theorem \ref{mosher}:}
Let $\Phi_1 \cdots \Phi_m$ be $m$ pseudo-anosov diffeomorphisms of
 $\Sigma$
 with
different sets of stable and unstable foliations. Let $H =
 \pi_1(\Sigma )$. Then there is an $n
 \geq 1$ such that the diffeomorphisms $\Phi_1^n, \cdots \Phi_m^n$
 generate a free group $F$ and the group $G$ given by the
exact sequence: \\
\begin{center}
$ 1 \rightarrow H \rightarrow G \rightarrow F \rightarrow 1$
\end{center}
\noindent is (strongly) hyperbolic relative to the
 maximal parabolic subgroups of the form $\ZZ \times F$. 

\medskip

We remark here that in Dahmani's combination theorem,
\cite{dahmani-comb}, an essential condition is
{\em acylindricity}. Again, in Alibegovic's combination theorem
\cite{alib-comb}, an essential assumption is the {\em compact
  intersection property}. Both acylindricity and the  compact
  intersection property prevent infinite (or even arbitrarily long
  chains of parabolics from occurring). This, to us, seemed
 a bit unsatisfactory,
  as the original motivation for the Bestvina-Feighn result came from
  Thurston's monster theorem (See \cite{kapovich-book}), and we wanted a
  generalization of the Bestvina-Feighn theorem that would cover the
  case of hyperbolic 3-manifolds with parabolics, particularly
  hyperbolic 3-manifolds of finite volume fibering over the circle. 
  The hypotheses in the present paper do allow for infinite chains of
  parabolics and covers the above case. Our 
emphasis here is geometric and so the main theorem  is stated in terms
of spaces rather than groups. 

\smallskip

{\bf Acknowledgements:} The authors are grateful to G. A. Swarup,
who was instrumental in bringing about this collaborative effort.
The work was completed during a visit of the second author
 to RKM Vivekananda University in February 2006.
We are also grateful to the referee for helpful comments and corrections.
 \footnote{After the submission of this paper, we learnt of the paper 
\cite{gautero-comb} by Gautero,
which gives a different proof of a result equivalent to Theorem 
\ref{strongcombingp}. 
A couple of points of difference between our work and 
\cite{gautero-comb}: We use Bestvina-Feighn's result \cite{BF}
directly, whereas an alternate proof of the combination
theorem of \cite{BF} is provided in 
\cite{gautero-comb}. However, in this paper, we 
provide in addition,  a converse (Theorem \ref{converse} ) to the main
combination theorem.}

\section{Relative Hyperbolicity}

In this section, we shall first recall certain notions of relative
hyperbolicity due to  
 Farb \cite{farb-relhyp} and Gromov \cite{gromov-hypgps}. 

\subsection{Electric Geometry}

Let $X$ be a path metric space. A collection of closed
 subsets $\HH = \{ H_\alpha\}$ of $X$ will be said to be {\bf uniformly
 separated} if there exists $\epsilon > 0$ such that
$d(H_1, H_2) \geq \epsilon$ for all distinct $H_1, H_2 \in \HH$.

\begin{definition} (Farb \cite{farb-relhyp})
The {\bf electric space} (or coned-off space) $\hhat{X}$
corresponding to the
pair $(X,\HH )$ is a metric space which consists of $X$ and a
collection of vertices $v_\alpha$ (one for each $H_\alpha \in \HH$)
such that each point of $H_\alpha$ is joined to (coned off at)
$v_\alpha$ by an edge of length $\half$. The sets $H_\alpha$ shall be
referred to as {\bf horosphere-like sets}. 
\label{el-space}
\end{definition}

A geodesic (resp. quasigeodesic) in $\hhat{X}$ will be referred to as
an electric geodesic (resp. quasigeodesic).

\begin{definition}
 A path $\gamma : I \rightarrow X$ in a path metric space $X$ is an
 {\bf ambient
K-quasigeodesic} if we have
\begin{center}
$L({\beta}) \leq K L(A) + K$
\end{center}
for any subsegment $\beta = \gamma |[a,b]$ and any rectifiable
path $A : [a,b] \rightarrow Y$ with the
same endpoints. (Here $L$ denotes length of path.)
\end{definition}

$N_R(Z)$ will denote the
$R$-neighborhood about the subset $Z$ in $X$.
 $N_R^e(Z)$ will denote the
$R$-neighborhood about the subset $Z$ in the electric metric.

Much of what Farb proved in \cite{farb-relhyp} goes through
under
considerably weaker assumptions than those of \cite{farb-relhyp}. In \cite{farb-relhyp} the theorems
were proven in the particular context of a pair $(X, {\mathcal{H}})$,
where $X$ is a Hadamard space of pinched negative curvature
 with the interiors of a family of
horoballs $\mathcal{H}$ removed. Then $\mathcal{H}$ can be regarded as
a collection of horospheres in $X$ separated by a minimum distance
from each other. In this situation, $X$ is not a hyperbolic metric
space itself, but is hyperbolic relative to a collection of separated
{\em horospheres}. Alternately let $H_h$ be the horoball corresponding
to the horosphere $H \in \mathcal{H}$. Let 
$X_h = X \bigcup_{H \in \mathcal{H}} H_h$ be the entire Hadamard manifold
of pinched negative curvature. Then the coned off space $\hhat{X}$
obtained by coning off the horospheres of $X$ is essentially 
equivalent to coned off space $\hhat{X_h}$
obtained by coning off the horoballs of $X_h$.

\smallskip

We consider
therefore a hyperbolic metric space $X$ and a collection $\mathcal{H}$
of
{\em (uniformly) $ C$-quasiconvex uniformly separated subsets}, i.e.
there exists $D > 0$ such that for $H_1, H_2 \in \mathcal{H}$, $d_X (H_1,
H_2) \geq D$. In this situation $X$ is weakly
hyperbolic relative to the
collection $\mathcal{H}$ in the sense that the coned off space 
$\hhat{X}$ is hyperbolic. The result in this form is due to Klarreich
\cite{klarreich}. However, the property of {\em Bounded Horosphere
Penetration (BHP)} or {\em Bounded Coset Penetration (BCP)} used by
Farb \cite{farb-relhyp} 
was not abstracted out
in Klarreich's proof as it was not necessary. What is essential for
BCP (or BHP) to go through has been abstracted out by Bowditch
\cite{bowditch-ct} \cite{bowditch-relhyp} in the case that the
collection $\mathcal{H}$ is a collection of geodesics or horocycles in
a Farey graph. (See also Bumagin \cite{bumagin}.)
But though these things are available at the level of
folklore, an explicit statement seems to be lacking.

The crucial
condition can be isolated as per the following definition \cite{brahma-ibdd}:

\smallskip

{\bf Definition:} A collection $\mathcal{H}$ of uniformly
$C$-quasiconvex sets in a $\delta$-hyperbolic metric space $X$
is said to be {\bf mutually D-cobounded} if 
 for all $H_i, H_j \in \mathcal{H}$, $\pi_i
(H_j)$ has diameter less than $D$, where $\pi_i$ denotes a nearest
point projection of $X$ onto $H_i$. A collection is {\bf mutually
  cobounded } if it is mutually D-cobounded for some $D$. 

\smallskip

{\em Mutual coboundedness} was proven for horoballs by Farb in Lemma 4.7 of
\cite{farb-relhyp} and by Bowditch in stating that the projection of
the link of a vertex onto another \cite{bowditch-relhyp} has bounded
diameter in the link. However, the comparability of 
intersection patterns in this context needs to be stated a bit more
carefully. We give the general version of Farb's theorem below and
refer to \cite{farb-relhyp} and Klarreich \cite{klarreich} for proofs.

\begin{lemma} (See Lemma 4.5 and Proposition 4.6 of
  \cite{farb-relhyp}, see also \cite{brahma-ibdd})
Given $\delta , C, D$ there exists $\Delta$ such that
if $X$ is a $\delta$-hyperbolic metric space with a collection
$\mathcal{H}$ of $C$-quasiconvex $D$-separated sets.
then,

\begin{enumerate}
\item {\it Electric quasi-geodesics electrically track hyperbolic
  geodesics:} Given $P > 0$, there exists $K > 0$ with the following
  property: Let $\beta$ be any electric $P$-quasigeodesic from $x$ to
  $y$, and let $\gamma$ be the hyperbolic geodesic from $x$ to $y$. 
Then $\beta \subset N_K^e ( \gamma )$. \\
\item $\gamma$ lies in a {\em hyperbolic} $K$-neighborhood of $N_0 ( \beta
  )$, where $N_0 ( \beta )$ denotes the zero neighborhood of $\beta$
  in the {\em electric metric}. \\
\item {\it Hyperbolicity:} The electric space
  $\hhat{X}$  is $\Delta$-hyperbolic. \\
\end{enumerate}
\label{farb1A}
\end{lemma}

Item (2) in the above Lemma is due to Klarreich \cite{klarreich},
where the proof is given for $\beta$ an electric geodesic, but the
same proof goes through for electric quasigeodesics without backtracking. 

The above Lemma motivates:

\begin{definition} \cite{farb-relhyp} \cite{bowditch-relhyp}
Let $X$ be a geodesic metric space and $\HH$ be a collection of
mutually disjoint uniformly separated subsets. Then $X$ is said to be
{\bf weakly hyperbolic} relative to the collection $\HH$, if the
electric space $\hhat{X}$ is hyperbolic.
\end{definition}

We shall need to give a general definition of geodesics and
quasigeodesics without backtracking. \\

{\bf Definitions:} Given a collection $\mathcal{H}$
of $C$-quasiconvex, $D$-separated sets and a number $\epsilon$ we
shall say that a geodesic (resp. quasigeodesic) $\gamma$ is a geodesic
(resp. quasigeodesic) {\bf without backtracking} with respect to
$\epsilon$ neighborhoods if $\gamma$ does not return to $N_\epsilon
(H)$ after leaving it, for any $H \in \mathcal{H}$. 
A geodesic (resp. quasigeodesic) $\gamma$ is a geodesic
(resp. quasigeodesic) {\bf without backtracking} if it is a geodesic
(resp. quasigeodesic) without backtracking with respect to
$\epsilon$ neighborhoods for some $\epsilon \geq 0$.

\smallskip

{\bf Note:} For the above lemma, the hypothesis is that $\mathcal{H}$
  consists of uniformly quasiconvex, mutually separated sets.
Mutual coboundedness has not yet been used. We introduce it in the next lemma.

\begin{lemma}  \cite{brahma-ibdd}
Suppose $X$ is a $\delta$-hyperbolic metric space with a collection
$\mathcal{H}$ of $C$-quasiconvex $K$-separated $D$-mutually cobounded
subsets. There exists $\epsilon_0 = \epsilon_0 (C, K, D, \delta )$ such that
the following holds:

Let $\beta$ 
  be an electric $P$-quasigeodesic without backtracking
and $\gamma$ a hyperbolic geodesic,
  both joining $x, y$. Then, given $\epsilon \geq \epsilon_0$
 there exists $D = D(P, \epsilon )$ such that \\
\begin{enumerate}
\item {\it Similar Intersection Patterns 1:}  if
  precisely one of $\{ \beta , \gamma \}$ meets an
  $\epsilon$-neighborhood $N_\epsilon (H_1)$
of an electrocuted quasiconvex set
  $H_1 \in \mathcal{H}$, then the length (measured in the intrinsic path-metric
  on  $N_\epsilon (H_1)$ ) from the entry point
  to the 
  exit point is at most $D$. \\
\item {\it Similar Intersection Patterns 2:}  if
 both $\{ \beta , \gamma \}$ meet some  $N_\epsilon (H_1)$
 then the length (measured in the intrinsic path-metric
  on  $N_\epsilon (H_1)$ ) from the entry point of
 $\beta$ to that of $\gamma$ is at most $D$; similarly for exit points. \\
\end{enumerate}
\label{farb2A}
\end{lemma}
 
Lemma \ref{farb2A} is essentially a paraphrasing of the BCP
property \cite{farb-relhyp} in terms of mutual coboundedness.
The above Lemma motivates the following definition:

\begin{definition} \cite{farb-relhyp} \cite{bowditch-relhyp}
Let $X$ be a geodesic metric space and $\HH$ be a collection of
mutually disjoint uniformly separated subsets such that
 $X$ is 
{\em weakly hyperbolic} relative to the collection $\HH$. If any pair of
electric quasigeodesics without backtracking starting and ending at
the same point have {\em similar intersection patterns} with {\em
  horosphere-like sets} (elements of $\HH$) then quasigeodesics
are said to satisfy {\bf Bounded Penetration} and 
 $X$ is said to be
{\em strongly hyperbolic} relative to the collection $\HH$.
\end{definition}

We summarize the two Lemmas \ref{farb1A} and \ref{farb2A}
as follows:

\noindent $\bullet$ If $X$ is a hyperbolic metric space and
$\mathcal{H}$ a collection of uniformly quasiconvex separated subsets,
then $X$ is hyperbolic relative to the collection $\mathcal{H}$. \\

\noindent $\bullet$ If $X$ is a hyperbolic metric space and
$\mathcal{H}$ a collection of uniformly quasiconvex mutually cobounded
separated subsets,
then $X$ is hyperbolic relative to the collection $\mathcal{H}$ and
satisfies {\em Bounded Penetration}, i.e. hyperbolic geodesics and
electric quasigeodesics have similar intersection patterns in the
sense of Lemma \ref{farb2A}.

\subsection{Partial Electrocution}

In this subsection, we indicate, following \cite{brahma-amalgeo},
 a modification of Farb's \cite{farb-relhyp}
 notion of
{\em strong relative hyperbolicity}
 and his construction of an electric metric, described earlier. The
 modification we shall discuss is called {\em partial electrocution} and will be used in proving the converse to the Strong Combination Theorem.
Most of this discussion is taken from \cite{brahma-amalgeo}.

\smallskip

We start with a few motivating examples:\\
\noindent 
{\bf Partial Electrocution} of a horosphere $H = \R^{n-1} \times \R$ 
will be defined as putting the
 zero metric in the $\R^{n-1}$ direction, and retaining the 
usual Euclidean metric in the other $\R$ direction.

In the partially electrocuted case, instead of coning all of a
horosphere down to a point 
we cone only horocyclic leaves of a foliation of the horosphere.
Effectively, therefore, we have a cone-line rather a cone-point. 

Let $Y$ be a convex simply connected hyperbolic n-manifold.
Let $\mathcal{B}$ denote a collection of horoballs. Let $X$ denote
$Y$ minus the interior of the horoballs in $\mathcal{B}$. Let 
$\mathcal{H}$ denote the collection of boundary horospheres.Then each
$H \in \mathcal{H}$ with the induced metric is isometric to a Euclidean
product $E^{n-2} \times L$ for an interval $L\subset \mathbb{R}$. 
Partially electrocute  each 
$H$ by giving it the product of the zero metric with the Euclidean metric,
i.e. on $E^{n-2}$ give the zero metric and on $L$ give the Euclidean
metric. The resulting space is exactly what one would get by gluing
to each $H$ the mapping cylinder of the projection of $H$ onto the $L$-factor.

\smallskip

This motivates the following scenario:

\begin{definition}
Let $(X, \HH , \GG , \LL )$ be an ordered quadruple such that the
following holds:

\begin{enumerate}
\item $X$ is (strongly) hyperbolic relative to a collection of subsets
$H_\alpha$, thought of as horospheres (and {\em not horoballs}). 
\item For each $H_\alpha$ there is a uniform large-scale
retraction $g_\alpha : H_\alpha \rightarrow L_\alpha$ to some
(uniformly) $\delta$-hyperbolic metric space $L_\alpha$, i.e. there
exist $\delta , K, \epsilon > 0$ such that for all $H_\alpha$ there exists
a $\delta$-hyperbolic $L_\alpha$ and a map 
$g_\alpha : H_\alpha \rightarrow L_\alpha$ with
$d_{L_\alpha} (g_\alpha (x), g_\alpha (y)) \leq Kd_{H_\alpha}(x,y)
+ \epsilon $ for all $x, y \in H_\alpha$. Further, we denote the
collection of such $g_\alpha$'s as $\GG$.  
\end{enumerate}
 The {\bf partially electrocuted space} or
{\em partially coned off space} corresponding to $(X, \HH , \GG , \LL)$ 
is 
obtained from $X$ by gluing in the (metric)
mapping cylinders for the maps
 $g_\alpha : H_\alpha \rightarrow L_\alpha$.
\end{definition}

In Farb's construction $L_\alpha$ is just a single point. However,
the notions and arguments of \cite{farb-relhyp} or Klarreich
 \cite{klarreich} 
go
through even in this setting. The metric, and geodesics and quasigeodesics
in the partially electrocuted space will be referred to as the 
partially electrocuted metric $d_{pel}$, and partially
electrocuted geodesics and quasigeodesics respectively. In this
situation, we conclude as in Lemma \ref{farb1A}:

\begin{lemma}
$(X,d_{pel})$ is a hyperbolic metric space and the sets $L_\alpha$
are uniformly quasiconvex.
\label{pel}
\end{lemma}

\noindent {\bf Note 1:} When $K_\alpha$ is a point, the last statement is a 
triviality.

\noindent {\bf Note 2:} $(X, d_{pel})$ is strongly hyperbolic relative to
the sets $\{ L_\alpha \}$. In fact the space obtained by electrocuting the
sets $L_\alpha$ in $(X,d_{pel})$ is just the space $(X,d_e)$ obtained
by (completely) 
electrocuting the sets $\{ H_\alpha \}$ in $X$.

\noindent {\bf Note 3:} The proof of Lemma \ref{pel} and other such
results below follow Farb's \cite{farb-relhyp} constructions. For
instance, consider a hyperbolic geodesic $\eta$ in a convex complete
simply connected  n-manifold $X$ with pinched negative curvature. 
Let $H_i$, $i = 1\cdots k$
be the partially electrocuted horoballs  it
meets. Let $N(\eta )$ denote the union of $\eta$ and $H_i$'s. Let $Y$
denote $X$ minus the interiors of the $H_i$'s. The
first step is to show that $N(\eta ) \cap Y$ is quasiconvex in $(Y,
d_{pel})$. To do this one takes a hyperbolic $R$-neighborhood of
$N(\eta )$ and projects
$(Y, d_{pel})$ onto it, using the hyperbolic projection. It was shown
by Farb in \cite{farb-relhyp} that the projections of all horoballs
 are uniformly
bounded in hyperbolic diameter. (This is essentially mutual
coboundedness). Hence, given $K$, choosing
 $R$  large enough, any path that goes out of
 an $R$-neighborhood of $N( \eta )$ cannot be a $K$-partially
 electrocuted
quasigeodesic. This is the one crucial step that allows the results of
\cite{farb-relhyp}, in particular, Lemma \ref{pel}
 to go through in the context of partially
electrocuted spaces.

\smallskip

As in Lemma \ref{farb2A}, partially electrocuted quasigeodesics
and geodesics without backtracking have the same intersection patterns
with {\em horospheres and boundaries of lifts of tubes} as electric geodesics
without backtracking.  
Further, since 
electric geodesics and hyperbolic quasigeodesics have similar intersection
 patterns with {\em horoballs and lifts of tubes} it follows that
 partially electrocuted 
quasigeodesics and hyperbolic quasigeodesics have similar intersection
patterns with {\em horospheres and boundaries of lifts of tubes}. We
state this formally below: 

\begin{lemma}
Given $K, \epsilon \geq 0$, there exists $C > 0$ such that the following
holds: \\
Let $\gamma_{pel}$ and $\gamma$ denote respectively a $(K, \epsilon )$
partially electrocuted quasigeodesic in $(X,d_{pel})$ and a hyperbolic
$(K, \epsilon )$-quasigeodesic in $(Y,d)$ joining $a, b$. Then $\gamma \cap X$
lies in a (hyperbolic) $C$-neighborhood of (any representative of) 
$\gamma_{pel}$. Further, outside of  a $C$-neighborhood of the horoballs
that $\gamma$ meets, $\gamma$ and $\gamma_{pel}$ track each other.
\label{pel-track}
\end{lemma}

\section{Trees of Hyperbolic Metric Spaces}

\subsection{Trees of Spaces: Hyperbolic and Relatively Hyperbolic}

We start with a notion closely related to one  introduced in \cite{BF}.

\medskip

\begin{definition}  A  tree (T) of hyperbolic (resp. strongly
  relatively  hyperbolic) metric spaces satisfying
the q(uasi) i(sometrically) embedded condition is a metric space $(X,d)$
admitting a map $P : X \rightarrow T$ onto a simplicial tree $T$, such
that there exist $\delta{,} \epsilon$ and $K > 0$ satisfying the following: \\
\begin{enumerate}
\item  For all vertices $v\in{T}$, 
$X_v = P^{-1}(v) \subset X$ with the induced path metric $d_v$ is a 
$\delta$-hyperbolic metric space (resp. a geodesic metric space $X_v$
strongly  hyperbolic relative to a collection $\HH_{v\alpha}$). Further, the
inclusions ${i_v}:{X_v}\rightarrow{X}$ 
are uniformly proper, i.e. for all $M > 0$, $v\in{T}$ and $x, y\in{X_v}$,
there exists $N > 0$ such that $d({i_v}(x),{i_v}(y)) \leq M$ implies
${d_v}(x,y) \leq N$.
\item Let $e$ be an edge of $T$ with initial and final vertices $v_1$ and
$v_2$ respectively.
Let $X_e$ be the pre-image under $P$ of the mid-point of  $e$.  
Then $X_e$ with the induced path metric is $\delta$-hyperbolic
(resp. a geodesic
 metric space $X_{e}$
strongly  hyperbolic relative to a collection $\HH_{e\alpha}$). 
\item There exist maps ${f_e}:{X_e}{\times}[0,1]\rightarrow{X}$, such that
$f_e{|}_{{X_e}{\times}(0,1)}$ is an isometry onto the pre-image of the
interior of $e$ equipped with the path metric.
\item ${f_e}|_{{X_e}{\times}\{{0}\}}$ and 
${f_e}|_{{X_e}{\times}\{{1}\}}$ are $(K,{\epsilon})$-quasi-isometric
embeddings into $X_{v_1}$ and $X_{v_2}$ respectively.
${f_e}|_{{X_e}{\times}\{{0}\}}$ and 
${f_e}|_{{X_e}{\times}\{{1}\}}$ will occasionally be referred to as
$f_{v_1}$ and $f_{v_2}$ respectively.
\item For a tree of strongly relatively hyperbolic spaces, we demand
  in addition, that the maps $f_{v_i}$ above ($i = 1, 2$) are {\bf
  strictly type-preserving}, i.e. $f_{v_i}^{-1}(H_{v_i \alpha})$, $i =
  1, 2$ (for
  any $H_{v_i \alpha} \in \HH_{v_i \alpha}$)
 is
  either empty or some $H_{e \alpha} \in \HH_{e \alpha}$.
\item For a tree of strongly relatively hyperbolic spaces, we demand
 that
 the
  coned off spaces are uniformly $\delta$-hyperbolic.
\end{enumerate}

$d_v$ and $d_e$ will denote path metrics on $X_v$ and $X_e$ respectively.
$i_v$, $i_e$ will denote inclusion of $X_v$, $X_e$ respectively into
$X$.

For a tree of relatively hyperbolic spaces, the sets $H_{v\alpha}$ and
$H_{e\alpha}$ will be referred to as {\bf horosphere-like vertex sets
  and edge sets} respectively.
\end{definition}

When $(X,d)$ is a tree (T) of strongly relatively hyperbolic metric
spaces, the {\em strictly type-preserving condition} (Condition 5
above)
ensures that we obtain an induced
  tree (T) (the same tree T) of {\em coned-off, or
 electric spaces}. We demand
  further \\

$\bullet$  {\bf qi-preserving electrocution condition: }
 the induced maps of the electric edge spaces into the
  electric vertex spaces $\hat{f_{v_i}} : \hhat{X_e} \rightarrow
  \hhat{X_{v_i}}$ ($i = 1, 2$) are uniform quasi-isometries.

The resulting tree of coned-off spaces will be called the {\bf induced
  tree of coned-off spaces}. The resulting space will be denoted as
  $\hhat{X}$. 

\noindent {\bf Definition:} A finite graph of (strongly) relatively
hyperbolic groups is said to satisfy Condition $C$, if the associated
tree of relatively hyperbolic Cayley graphs satisfies Condition
$C$. Here $C$ will be one of the following: \\
\begin{enumerate}
\item the qi-embedded condition 
\item  the strictly type-preserving
condition
\item the qi-preserving electrocution condition
\item the induced tree of coned-off spaces
satisfies the {\bf hallways flare} condition (See below)
\item the {\bf cone-bounded hallways strictly flare} condition (See below)
\end{enumerate}

{\bf Remark:}  Strictly speaking, this induced tree
  exists
for any collection of vertex and edge spaces satisfying the {\em
  strictly type-preserving condition}. Hyperbolicity is not essential
  for the existence of the induced tree of spaces.

The {\bf cone locus} of $\hhat{X}$,
the induced tree (T) of coned-off spaces, is
the graph (in fact a forest) whose vertex set $\VV$ consists of the
cone-points in the vertex set of $\hhat{X}$
 and whose edge-set $\EE$
consists of the
cone-points in the edge set
of $\hhat{X}$. The incidence relations are dictated by
 the incidence relations in $T$.

Note that connected components of the cone-locus can be naturally
identified with sub-trees of $T$. Each such connected component of the
cone-locus will be called a {\bf maximal cone-subtree}. The collection
of {\em maximal cone-subtrees} will be denoted by $\TT$ and elements
of $\TT$ will be denoted as $T_\alpha$. Further, each maximal
cone-subtree $T_\alpha$ naturally gives rise to a tree $T_\alpha$ of
horosphere-like subsets depending on which cone-points arise as
vertices and edges of $T_\alpha$. The metric space that $T_\alpha$
gives rise to will be denoted as $C_\alpha$ and will be referred to as
a {\bf maximal cone-subtree of horosphere-like spaces}. The collection
of $C_\alpha$'s will be denoted as $\CC$. \\
{\bf Note:} Each $T_\alpha$ thus appears in two guises:\\
\begin{enumerate}
\item as a subset of $\hhat{X}$
\item as the underlying tree of $C_\alpha$
\end{enumerate}
We shall have need for both these interpretations.

\subsection{The Bestvina-Feighn Flare Condition}

Next, we would like to recall the essential condition (due to
 Bestvina and Feighn
\cite{BF}) ensuring
 hyperbolicity of a tree of spaces. We retain the terminology.

\begin{definition}  
A disk $f : [-m,m]{\times}{I} \rightarrow 
X$ is a hallway of length $2m$ if it satisfies:

\begin{enumerate}
\item $f^{-1} ({\cup}{X_v} : v \in T) = \{-m,  \cdots , m \}{\times}
I$
\item $f$ maps $i{\times}I$ to a geodesic in  $X_v$ for some vertex
space.
\item $f$ is transverse, relative to condition (1) to $\cup_e X_e$.

\end{enumerate}
\end{definition}

\begin{definition} A hallway is $\rho$-thin if 
$d({f(i,t)},{f({i+1},t)}) \leq \rho$ for all $i, t$.

A hallway is {\em $\lambda$-hyperbolic} if \\

$\lambda l(f(\{ 0 \} \times I)) \leq $ max $\{ l(f(\{ -m \} \times I)),
l(f(\{ m \} \times I))$

A hallway is {\em essential} if the edge path in $T$ 
resulting from projecting $X$
 onto $T$ does not backtrack (and is therefore a geodesic segment in
 the tree $T$).

An essential  hallway of length $2m$ is {\bf cone-bounded} if
$f(i \times {\partial I})$
lies in the cone-locus for $i = \{ -m, \cdots , m\}$.
\end{definition}

\begin{definition} {\bf Hallways flare condition:}
The tree of spaces, $X$, is said to satisfy the {\em hallways flare}
condition if there are numbers $\lambda > 1$ and $m \geq 1$ such that
for all $\rho$ there is a constant $H(\rho )$ such that  any
$\rho$-thin essential hallway of length $2m$ and girth at least $H$ is
$\lambda$-hyperbolic.
\end{definition}

\begin{definition} {\bf Cone-bounded hallways strictly flare
 condition:} 
The tree of spaces, $X$, is said to satisfy the {\em hallways flare}
condition if there are numbers $\lambda > 1$ and $m \geq 1$ such that
 any
cone-bounded hallway of length $2m$  is
$\lambda$-hyperbolic.
\end{definition}

The main theorem of Bestvina and Feighn follows
(though this is stated in \cite{BF} for groups, the proof nowhere
requires uniform properness of the
space). Bowditch \cite{bowditch-ct} notes the equivalence of the hallways flare condition with the hyperbolicity of the tree of spaces. 

\begin{theorem} \cite{BF} pp 85-86 \cite{bowditch-ct}
Let $X$ be a tree of hyperbolic metric spaces satisfying the {\em
  q.i.-embedded condition} and the {\em hallways flare
  condition}. Then $X$ is hyperbolic. \\
Conversely, if $X$ is hyperbolic, then hallways flare.
\label{BF}
\end{theorem}

Apart, from Theorem \ref{BF} above, we shall need one more simple
 observation.

\begin{lemma}
Suppose that $\hhat{X}$ is hyperbolic. Then
maximal cone-subtrees $T_\alpha$ are uniformly
quasi-convex in $\hhat{X}$.
\label{subtree-qc}
\end{lemma}

\noindent {\bf Proof:} Let $P: X \rightarrow T$ be the natural
projection map of the tree of spaces to the underlying sub-tree. Then
$P$ induces $P': \hhat{X} \rightarrow T$ as $\hhat{X}$ may be regarded
as (the same) tree (T) of coned-off spaces. $P'$ is distance
non-increasing. Further, restricted to each $T_\alpha$, $P'$ is an
isometry. 

Also note that any path from $x \in H_{v1}$ to $y \in H_{v2}$ in
$\hhat{X}$ has length not less than $d_T (P'(x),P'(y))$, where $d_T$
is the natural metric on $T$ and $v_1, v_2 \in T_\alpha$.

Now suppose that $x, y \in T_\alpha \subset \hhat{X}$.  Let $\gamma
\subset T_\alpha \hhat{X}$ be the geodesic in $T_\alpha$ joining $x, y$.
It therefore follows that for any $x, y \in T_\alpha \subset \hhat{X}$
and any path $A$ joining $x, y \in \hhat{X}$, 

\begin{center}
$l(A) \leq d_T(P'(x), P'(y)) = l(\gamma)$
\end{center}

Hence $\gamma$ is quasi-isometrically (in fact isometrically) embedded
in $\hhat{X}$ and hence a geodesic in $\hhat{X}$. The lemma
follows. $\Box$

\section{The Combination Theorem}

\subsection{Weak Combination Theorem}

We start with the following

\begin{theorem} {\bf Weak Combination Theorem}
Let $X$ be a tree ($T$) of strongly relatively hyperbolic spaces
satisfying 
\begin{enumerate}
\item the qi-embedded condition 
\item  the strictly type-preserving
condition
\item the qi-preserving electrocution condition
\item the induced tree of coned-off spaces
satisfies the {\em hallways flare} condition
\end{enumerate}
Then $X$  is {\em
  weakly hyperbolic} relative to the family $\CC$ of {\em maximal
  cone-subtrees of horosphere-like spaces}.
\label{weakcombin}
\end{theorem}

\noindent {\bf Proof:} As usual let $\hhat{X}$ denote the induced tree
($T$) of coned-off spaces, $\TT$ denote the family of maximal
cone-subtrees $T_\alpha \subset \hhat{X}$. Let $\hhat{\hhat{X}}$
denote $\hhat{X}$ with the family of sets $\TT$ coned off (i.e. 
vertices $v_\alpha$
are introduced, one each for each $T_\alpha$, and joined to points of
the corresponding $T_\alpha$ by edges of length $\half$.)

Since vertex and edge-spaces are strongly relatively hyperbolic,
then by item (6) in the
definition of a tree of strongly relatively hyperbolic spaces,
 $\hhat{X}$ is a tree of (uniformly) hyperbolic metric spaces.

 By the
{\em qi-preserving electrocution condition}, the induced tree of
coned-off spaces satisfies the
qi-embedded condition. 

By the {\em hallways flare condition} and Theorem \ref{BF}, $\hhat{X}$
is a hyperbolic metric space. 

By Lemma \ref{subtree-qc}, the sets $T_\alpha \in \TT$ are uniformly
quasiconvex and uniformly separated. 

Hence by Lemma \ref{farb1A}, $\hhat{X}$ is weakly
hyperbolic relative to the sets $T_\alpha \in \TT$, i.e.
 $\hhat{\hhat{X}}$ is a hyperbolic metric space.

Let $\hhat{X_1}$ denote the space obtained from $X$ by coning off
{\em maximal cone-subtrees of horosphere-like sets}. Then $\hhat{X_1}$
is  quasi-isometric to  $\hhat{\hhat{X}}$. To see this, one notes that
 $\hhat{\hhat{X}}$ is obtained from $X$ by first coning-off (or partially
electrocuting) $C_\alpha$'s, 
the maximal cone-subtrees of horosphere-like sets, to
maximal cone-subtrees $T_\alpha$. This gives rise to
$\hat{X}$. Further coning off the $T_\alpha$'s gives
$\hhat{\hhat{X}}$. On the other hand, $\hhat{X_1}$ is obtained from
$X$ by coning off (or  completely electrocuting)
the $C_\alpha$'s to points  in one step. The two constructions clearly
give quasi-isometric spaces.

Hence $\hhat{X_1}$ is hyperbolic, i.e. $X$ is weakly hyperbolic
relative to the collection of sets $C_\alpha \in \CC$. $\Box$

\subsection{Strong Combination Theorem}

Under the additional  {\em cone-bounded hallways strictly
flare condition}, we would now like to prove a stronger version of the
combination Theorem \ref{weakcombin}, i.e. $X$ is {\bf strongly
  hyperbolic} relative to the collection of $C_\alpha \in \CC$:

By Lemma \ref{farb2A} and Theorem \ref{weakcombin}, it suffices to
show that the sets $T_\alpha \subset \hhat{X}$ are mutually
co-bounded.
Most of the rest of this subsection is devoted to proving {\em mutual
  coboundedness}. 

The next lemma follows easily from stability of quasigeodesics \cite{GhH}
\cite{gromov-hypgps}. (See, for instance Lemma 4.1.1 of \cite{mitra-thesis}.)

\begin{lemma}
Give $\delta , C$, there exist $D, K, \epsilon$ such that the
following holds: \\
Let $(X,d)$ 
be a $\delta$-hyperbolic metric space and $Y$ a
$C$-quasiconvex subset. Let $\pi$ be a nearest-point retraction of $X$
onto $Y$. Let $x, y \in X$ such that $d(\pi (x),\pi (y)) \geq D$. Then
$[x,\pi (x)]\cup [\pi (x), \pi (y)]\cup [\pi (y), y]$ is a $(K,
\epsilon )$-quasigeodesic.
\label{qg1}
\end{lemma}

We use Lemma \ref{qg1} below:

\begin{cor}
Given $\delta$, $C$, there  exist $D, M$ such that the following holds:\\
Suppose that  $Y, Z$ are
$C$-quasiconvex subsets of  a $\delta$-hyperbolic metric space
$(X,d)$. Let $\pi$ denote nearest point projection onto $Y$. If $\pi
(Z)$ 
has diameter greater than $D$, then $\pi (Z)$ lies in an
$M$-neighborhood of $Z$. 
\label{largeprojn}
\end{cor}

\noindent {\bf Proof:} Let $x, y \in Z$. By Lemma \ref{qg1}, there
exist
$D_0, K, \epsilon$ (depending on $\delta$, $C$) such that if
$d(\pi (x),\pi (y)) \geq D_0$, then
$[x,\pi (x)]\cup [\pi (x), \pi (y)]\cup [\pi (y), y] = \gamma$ is a $(K,
\epsilon )$-quasigeodesic. Since $Z$ is $C$-quasiconvex, $\gamma$ lies
in an $M_0 = M_0(K, \epsilon , C, \delta )$-neighborhood of $Z$. 

Now, choose $a, b$ in $Y, Z$ respectively, such that $d(a,b) =
d(Y,Z)$. 
From the previous paragraph, we deduce that if $z \in Z$ such that
$d(\pi (z), a) \geq D_0$, then $\pi(z)$ lies in an $M_0$-neighborhood of
$Z$. Taking $D  = 2D_0$ and $M = M_0 + D_0$, we are through. $\Box$

The following Proposition deduces mutual coboundedness of 
maximal cone-subtrees
$T_\alpha$ in $\hhat{X}$ (obtained by partially electrocuting 
 maximal cone-subtrees of horosphere-like sets) from hyperbolicity
of $\hhat{X}$ (established for instance in  Theorem \ref{weakcombin} )
and {\bf cone-bounded hallways strictly flare} condition.
This will be used in proving the Strong Combination Theorem.  

\begin{prop}
Suppose that the tree of coned-off spaces $\hhat{X}$ is hyperbolic and
that the {\bf cone-bounded hallways strictly flare} condition is
satisfied.
Then
there exists $D \geq 0$ such that the family of maximal cone-subtrees
$T_\alpha$ in $\hhat{X}$ is $D$-cobounded. 
\label{con-tree-cobdd}
\end{prop}

\noindent {\bf Proof:} Suppose not. Then, by Corollary
\ref{largeprojn}, there exists $M \geq 0$ such that
for any $D \geq 0$,
there exist maximal cone-subtrees $T_1, T_2$ and a connected subtree $T_3
\subset T_1$ such that $T_3$ has diameter  greater than $D$ and lies
in an $M$-neighborhood of $T_2$. Hence, there exist geodesic
edge-paths $\gamma_1 \subset T_1$ and $\gamma_2 \subset T_2$ such that
each $\gamma_i$ has length greater than $(D-2M)$ and such that they
lie in an
$M$-neighborhood of each other in $\hhat{X}$. 

Next, $P(\gamma_i) \subset T$
(where $T$ is the tree underlying $\hhat{X}$, 
the tree (T) of coned-off spaces and $P: \hhat{X} \rightarrow T$
is the natural projection). Also,  $P(\gamma_i)$ ($i = 1, 2$) is
abstractly isomorphic to $\gamma_i$ as an edge-path, since $P :
\hhat{X} \rightarrow T$ is an isometry restricted to each $T_i$. 
Then $P (\gamma_i)$ may be regarded as geodesic paths in $T$
having lengths greater than $(D-2M)$ and lying in an $M$-neighborhood
of each other (since $P$ does not increase distances). This means that
$P (\gamma_i)$ (for $i= 1, 2$)
must overlap over an interval of length at least $(D-4M)$.

Let $\alpha_i \subset \gamma_i$ be paths having length at least
$(D-4M)$ with $P(\alpha_1) = P(\alpha_2)$. Then there exists $M_1 =
M_1(M)$
and a
{\em cone-bounded
 hallway}
$\Delta : [-m,m] \times I \rightarrow \hhat{X}$ with $2m \geq (D-4M)$
such that \\
\begin{enumerate}
\item $\Delta ([-m,m] \times \{ 0 \} = \alpha_1$ 
\item $\Delta ([-m,m] \times \{ 1 \} = \alpha_2$
\item each $\Delta (j \times I)$ has length less than $M_1$
\end{enumerate} 

Since $D$, and hence $m$  can be arbitrarily large, while $M$ (and
hence $M_1$) are fixed,it follows that for any given $\lambda > 1$,
there exists a hallway $\Delta$, which is not $\lambda$-hyperbolic.
This violates the {\em cone-bounded hallways strictly flare}
condition.
Hence, by contradiction, 
there exists $D \geq 0$ such that the family of maximal cone-subtrees
$T_\alpha$ in $\hhat{X}$ is $D$-cobounded. 
$\Box$

We are now in a position to prove:\\

\begin{theorem} {\bf Strong Combination Theorem}
Let $X$ be a tree ($T$) of strongly relatively hyperbolic spaces
satisfying 
\begin{enumerate}
\item the qi-embedded condition 
\item  the strictly type-preserving
condition
\item the qi-preserving electrocution condition
\item the induced tree of coned-off spaces
satisfies the {\bf hallways flare} condition
\item the {\bf cone-bounded hallways strictly flare} condition.
\end{enumerate}
Then $X$  is {\em
  strongly hyperbolic} relative to the family $\CC$ of {\em maximal
  cone-subtrees of horosphere-like spaces}.
\label{strongcombin}
\end{theorem}

\noindent {\bf Proof:} By Theorem \ref{weakcombin}, we know that $X$
is weakly hyperbolic relative to the family $\CC$ of {\em maximal
  cone-subtrees of horosphere-like spaces}.

This is equivalent to saying that $\hhat{X}$ is weakly hyperbolic
  relative to the family $\TT$ of {\em maximal 
  cone-subtrees} $T_\alpha \subset \hhat{X}$.

By the  {\bf cone-bounded hallways strictly flare} condition and 
Proposition \ref{con-tree-cobdd}, we see that the family $\TT$ is
mutually cobounded.

Hence by Lemma \ref{farb2A}, we conclude that 
 $\hhat{X}$ is strongly hyperbolic
  relative to the family $\TT$ of {\em maximal 
  cone-subtrees} $T_\alpha \subset \hhat{X}$. Equivalently,
 $X$
is strongly hyperbolic relative to the family $\CC$ of {\em maximal
  cone-subtrees of horosphere-like spaces}. $\Box$ 

\medskip

Recall that a finite graph of (strongly) relatively
hyperbolic groups is said to satisfy  a Condition $C$, if the associated
tree of relatively hyperbolic Cayley graphs also satisfies Condition
$C$. The resulting group will be denoted as $G$. 
A quotient of {\em  maximal
  cone-subtrees of horosphere-like spaces} in this case, is called
a  {\em maximal
  cone-subgraph of horosphere-like subgroups}. Note that such a
subgraph gives rise to a subgroup of $G$. We shall refer to such
subgroups as  {\bf maximal
  parabolic subgroups}. Recall that we are using the convention that
all parabolic subgroups are of infinite index.
As an immediate consequence of Theorem \ref{strongcombin}, we have the
following:

\begin{theorem} {\bf Strong Combination Theorem for Graphs of Groups}
Let $G$ be a finite graph ($\Gamma$) of strongly relatively hyperbolic groups
satisfying 
\begin{enumerate}
\item the qi-embedded condition 
\item  the strictly type-preserving
condition
\item the qi-preserving electrocution condition
\item the induced tree of coned-off spaces
satisfies the {\bf hallways flare} condition
\item the {\bf cone-bounded hallways strictly flare} condition.
\end{enumerate}
Then $G$  is {\em
  strongly hyperbolic} relative to the family $\CC$ of {\em maximal
  parabolic subgroups}.
\label{strongcombingp}
\end{theorem}

\subsection{Converse to the Strong Combination Theorem}

 Recall that the  {\bf partially electrocuted space} or
{\em partially coned off space} corresponding to a quadruple
$(X, \HH , \GG , \LL)$ 
is obtained from $X$ by gluing in the (metric)
mapping cylinders for the maps  $g_\alpha : H_\alpha \rightarrow L_\alpha$.
Note that from Theorem 
\ref{strongcombin}, it follows that $\hhat{X}$ is obtained from $X$ by
partially electrocuting each $C_\alpha$.
Here 
\begin{enumerate}
\item $H_\alpha = C_\alpha$ and $\HH = \CC$
\item $L_\alpha = T_\alpha$  and $\LL = \TT$
\item $g_\alpha : C_\alpha
 \rightarrow T_\alpha $ collapses  $C_\alpha$, the tree
of  horosphere-like spaces to the underlying tree
 $T_\alpha$. 
\end{enumerate}

\begin{theorem} {\bf Converse to Strong Combination Theorem }
Let $X$ be a tree ($T$) of strongly relatively hyperbolic spaces
satisfying 
\begin{enumerate}
\item the qi-embedded condition 
\item  the strictly type-preserving
condition
\item the qi-preserving electrocution condition
\item $X$  is {\em
  strongly hyperbolic} relative to the family $\CC$ of {\em maximal cone-subtrees} of horosphere-like spaces
\end{enumerate}
Then
the induced tree of coned-off spaces
satisfies the {\bf hallways flare} condition
and the {\bf cone-bounded hallways strictly flare} condition.
\label{converse}
\end{theorem}

\noindent {\bf Proof:} As usual let $C_\alpha$ denote maximal cone sub-trees ($T_\alpha$) of horosphere-like sets.
By Lemma \ref{pel}, the induced tree of coned-off spaces
$\hhat{X}$, obtained by partially electrocuting each
$C_\alpha$ to $T_\alpha$. Then by the converse part of Theorem \ref{BF}, hallways, including cone-bounded hallways
flare.

It remains to show that cone-bounded hallways {\it strictly}
flare. Suppose not. Then there exists $D_0$ such that for all $N \in \mathbb{N}$, there exist cone-bounded hallways of length greater than $N$, bounded by "vertical" (parametrized)
geodesics
$\lambda_1, \lambda_2$ in distinct cone-subtrees $T_1, T_2$
respectively such that $d(\lambda_1(i), \lambda_2(i)) \leq D_0$ for all $i = 0, \cdots N$. Let $\mu_0, \mu_N$ denote
"horizontal" paths in the hallway joining $\lambda_1(i), \lambda_2(i)$ for $i= 0, N$. Hence, there exist points
$a_j, b_j$ ($j=0,N$) lying on the
corresponding cones $C_1\cap \mu_j, C_2\cap\mu_j$, 
respectively such that $d(a_j, b_j) \leq D_0$. Then
we have two paths:\\
$\bullet$ $\sigma_1$ starts at $a_0$, moves to $\lambda_1(0) $
(by a cone-edge of length $\frac{1}{2}$), proceeds to 
$\lambda_1(N) $ and exits to $a_N$
(again by a cone-edge of length $\frac{1}{2}$) \\
$\bullet$ $\sigma_2$ starts at $a_0$, moves to $b_0$
by a path of length $\leq D_0$, then to $\lambda_2(0) $
(by a cone-edge of length $\frac{1}{2}$), proceeds to 
$\lambda_2(N) $, exits to $b_N$
(again by a cone-edge of length $\frac{1}{2}$) and then goes
to $a_N$ by a path of length $\leq D_0$. \\
$\sigma_1$ has length $1$ in the electric metric on $\hhat{\hhat{X}}$ and $\sigma_2$ has length $\leq 2D_0 + 1$. Since $D_0$ is fixed, both $\sigma_1$ and $\sigma_2$ are uniform
quasigeodesics beginning and ending at the same point, but have
manifestly different intersection patterns with $C_1, C_2$.
This contradicts Lemma \ref{farb2A} and hence $X$ cannot be
strongly hyperbolic relative to the collection $C_\alpha$.
This final contradiction proves the theorem. $\Box$

\subsection{Examples}

The first class of examples are hyperbolic 3-manifolds fibering over
the circle with fiber a punctured hyperbolic surface $\Sigma$.
All the conditions of Theorem \ref{strongcombingp}
are satisfied in this case of a
3-manifold fibering over the circle with fiber a punctured surface.
The one condition that needs checking is the {\em hallways flare
  condition} for the induced tree (in fact line) of coned-off spaces.
This fact is due to  Bowditch (see \cite{bowditch-ct} Section 6). We
give here a somewhat different argument based on work of Mosher
\cite{mosher-hbh}. 

In \cite{mosher-hbh}, Mosher constructs
  examples of exact sequences of hyperbolic groups of the form: \\
\begin{center}
$ 1 \rightarrow H \rightarrow G \rightarrow F \rightarrow 1$
\end{center}
where $H$ is a closed surface group, $G$ is hyperbolic
 and $F$ is free. (This construction
is modified by Bestvina, Feighn and Handel \cite{BFH-lam} to the case
 where $H$ is a free group.)

We shall modify Mosher's argument slightly to make it work for
punctured surfaces. 

Let $\Phi$ be a pseudoanosov diffeomorphism
of  a punctured hyperbolic surface $\Sigma$. Taking a suitable power
of
$\Phi$ if necessary, we may assume that $\Phi$ fixes all the
punctures. The stable and unstable foliations of $\Phi$ give rise to a
piecewise Euclidean metric on $\Sigma$. This metric is incomplete at
the punctures. Complete it to get a surface with boundary $\Sigma_B$,
which may be thought of as the blow-up of $\Sigma$ at the
punctures. Equip $\Sigma_B$ with a pseudo-metric which is zero on all
the boundary components and equal to the piecewise Euclidean metric
elsewhere. This metric is discontinuous at the boundary, but this is
not important. This is essentially the {\em electric metric} on
$\Sigma_B$.

Then any electric geodesic $\lambda$ 
in $\Sigma_B$ is the union of two types of segments: \\

\begin{enumerate}

\item Geodesics in the piecewise Euclidean metric meeting the boundary
  at right angles. Let $\lambda_{eu}$ denote this union.
\item segments lying along the boundary.

\end{enumerate}

The total length of such an electric geodesic is the sum of the
lengths of the Euclidean pieces. The projection of the union of the
Euclidean pieces $\lambda_{eu}$
 onto the stable and unstable foliations will be denoted by
 $\lambda_{eus}$ and $\lambda_{euu}$ respectively. Abusing notation
 slightly, we assume that $\lambda_{eu}$, $\lambda_{eus}$,
 $\lambda_{euu}$ denote the respective lengths also. Then 
$max (\lambda_{eus},
 \lambda_{euu}) \geq \frac{1}{2} \lambda_{eu}$. Let $\Phi
 (\lambda_{eu})$
denote the image of $\lambda_{eu}$ under $\Phi$. If we assume that the
stable and unstable foliations meet the boundary components of
$\Sigma_B$ at right angles, then it can be easily shown that for any
given $k > 1$, there is
an $n$ (depending on the stretch factor of $\Phi$) such that \\
$max (\Phi(\lambda_{eus}),
 \Phi^{-1}(\lambda_{euu})) \geq k \lambda_{eu}$ and hence, \\
$max (\Phi(\lambda_{eu}),
 \Phi^{-1}(\lambda_{eu})) \geq k \lambda_{eu}$.

This proves the  one condition that needed checking, viz.
 the {\em hallways flare
  condition} for the induced tree (in fact line) of coned-off spaces.

More generally, we may take any $m$ (equal to two below for
 concreteness)
 pseudoanosov diffeomorphisms
 $\Phi, \Psi$ with different stable and unstable foliations. Then
 generalizing the above construction, we can prove the following
 generalization of  an essential Lemma of Mosher:{\em 3 out of 4
 stretch}. 

\begin{lemma} For any $k > 1$, there exists $n > 0$ such that for any
 electric geodesic $\lambda$ in $\Sigma_B$, at least three of the four
 elements $\Phi^n, \Phi^{-n}, \Psi^n, \Psi^{-n}$ stretch $\lambda$ by
 a factor of $\lambda$. 
\label{3of4}
\end{lemma}

 Thus we get an
 exact sequence of groups of the form: \\
\begin{center}
$ 1 \rightarrow H \rightarrow G \rightarrow F \rightarrow 1$
\end{center}
where $H$ is a punctured surface group, and 
 and $F$ is free.

The Cayley graph of $G$ may thus be regarded as a tree $T$
of hyperbolic spaces, where $T$ arises as the Cayley graph of
$F$. Also, the maximal parabolic subgroups
 here correspond exactly to the peripheral
subgroups, i.e. the cusp groups. 
Maximal cone-subtrees are each isometric to $T$. Then a tree
$T$ of maximal parabolics would correspond to
$\mathbb{Z} \times F$.
Lemma \ref{3of4}
shows that the induced tree of coned off spaces is hyperbolic.

Thus, we obtain from Theorem \ref{strongcombingp}

\begin{theorem}
Let $\Phi_1 \cdots \Phi_m$ be $m$ (orientation-preserving)
pseudo-anosov diffeomorphisms of
 $\Sigma$
 with
different sets of stable and unstable foliations. Let $H =
 \pi_1(\Sigma )$. Then there is an $n
 \geq 1$ such that the diffeomorphisms $\Phi_1^n, \cdots \Phi_m^n$
 generate a free group $F$ and the group $G$ given by the
exact sequence: \\
\begin{center}
$ 1 \rightarrow H \rightarrow G \rightarrow F \rightarrow 1$
\end{center}
is (strongly) hyperbolic relative to the
 maximal parabolic subgroups of the form 
$\mathbb{Z} \times F$. 
\label{mosher}
\end{theorem}

For $n = 1$, we get back Bowditch's theorem \cite{bowditch-ct}.

\subsection{Applications, Consequences and Problems}
Theorem \ref{strongcombin} and Theorem \ref{strongcombingp} open
 up the possibility of generalizing
several theorems about hyperbolic groups to (strongly) relatively
 hyperbolic groups. \\ 
 {\bf 1) Cannon-Thurston Maps}: In \cite{mitra-trees}, the
first author
 proved the existence of Cannon-Thurston maps for trees of hyperbolic
 metric spaces.  In \cite{brahma-pared}, he  generalized this theorem
  to the  relatively hyperbolic case
 under the additional assumption that the tree of
  spaces gives rise to a hyperbolic 3-manifold of bounded geometry
  whose core is incompressible away from cusps. In \cite{mahan-pal}, 
Mj-Pal prove the existence of Cannon-Thurston maps for the situation
 discussed in this paper, viz. trees of (strongly) relatively hyperbolic
 trees of metric spaces that are (strongly) relatively hyperbolic. 
\\ 
{\bf 2) Strongly relatively hyperbolic Extensions of Groups:} A
 Theorem of Mosher \cite{mosher-hypextns} says that if an exact sequence
 of groups of the form \\
\begin{center}
$ 1 \rightarrow H \rightarrow G \rightarrow Q \rightarrow 1$
\end{center}
exists, where $H$ is hyperbolic, then there exists a quasi-isometric
section of $\Gamma_Q$ into $\Gamma_G$ exists. In particular, if $G$ is
hyperbolic, so is $Q$. The essential technique is to use the action of
$Q$ on the boundary $\partial H$ of $H$.A fact (due to Gromov
\cite{gromov-hypgps}) that is used is that the space of triples of
points on the boundary of a hyperbolic group $H$ is
quasi-isometric to $ \Gamma_H$. An analogous result is shown by Pal in \cite{pal-thesis}.
\\ 
 {\bf 3) Heights of Groups:} In \cite{GMRS}, Gitik, Mitra, Rips and
 Sageev show that quasiconvex subgroups of hyperbolic groups have
 finite height and finite width. A partial converse was obtained by
 the first author in \cite{mitra-ht} for groups splitting over
 subgroups. This converse was used by Swarup in \cite{swarup-weakhyp}
 to prove a weak hyperbolization theorem. All three theorems should
 have analogues in the (strongly) relatively hyperbolic world. Hruska
 and Wise \cite{hruska-wise} have already shown the finiteness of
 height and width of quasiconvex subgroups of relatively hyperbolic groups.

\bibliography{relhyp}
\bibliographystyle{alpha}

\end{document}